\documentclass[preprint,3p,12pt]{elsarticle}

\usepackage{color}
\usepackage{amssymb}
\usepackage{amsmath}
\usepackage{amsthm}

\usepackage{epsfig}
\usepackage{ulem}
\usepackage{epic}
\usepackage{tabls}
\usepackage{booktabs}
\usepackage{threeparttable}
\usepackage{pict2e}

\usepackage{float}
 \usepackage{subfig}

\usepackage{hyperref}

\hypersetup{backref,colorlinks=true,linkcolor=blue}

\usepackage[vlined]{algorithm2e}

\begin{document}
\begin{frontmatter}



\title{Multilevel Iteration Method     for    Binary Stochastic Transport Problems}

\author[ncsu,ncsu1]{Dmitriy Y. Anistratov}
\address[ncsu]{Department of Nuclear Engineering,
North Carolina State University Raleigh, NC}
\address[ncsu1]{anistratov@ncsu.edu}

\begin{abstract}
This paper presents  an iteration method for solving linear particle transport problems
in binary stochastic mixtures. It is based on nonlinear projection approach. The method is defined by a hierarchy of equations consisting of
the high-order transport equation for materials,
low-order Yvon-Mertens equations for conditional ensemble average  of the material partial scalar fluxes,
 and low-order quasidiffusion equations for
  the  ensemble average  of the  scalar flux and current.
  The multilevel system of equations is solved by means of an iterative algorithm with the $V$-cycle. The  iteration method  is analyzed on a set of numerical test problems.
 \end{abstract}


\end{frontmatter}

\section{Introduction}

We consider the Levermore-Pomraning (LP) particle transport model
 for binary stochastic mixtures (BSM) in 1D slab geometry
 \cite{Levermore-1986,mla-jqsrt-1989,pomraning-1990}.
The materials are randomly distributed as alternating layers
 with  homogeneous Markov mixing statistics.
The mean width of  layers of  material $\ell$  is equal to~$\lambda_{\ell}$.
The BSM  transport equation with the LP closure in 1D slab geometry is given by
\begin{multline} \label{bsm-te}
\mu \frac{\partial}{\partial  x} \Big( p_{\ell} \psi_{\ell}  (x, \mu) \Big)
+ \sigma_{t, \ell} p_{\ell} \psi_{\ell}  (x, \mu) =
\frac{1}{2}\sigma_{s, \ell} p_{\ell} \int_{-1}^{1} \psi_{\ell}  (x, \mu') d\mu' \\
+ \frac{|\mu|}{\lambda_{\ell'}}p_{\ell'} \psi_{\ell'}  (x, \mu)
-   \frac{|\mu|}{\lambda_{\ell}}p_{\ell} \psi_{\ell}  (x, \mu)
+ \frac{1}{2}p_{\ell} q_{\ell}(x)\, ,
\end{multline}
\[
 \ell=   1,2 \, , \quad \ell' \ne \ell \, , \quad   x \in [0,X]  \, , \quad  \mu \in [-1,1] \,
\]
 with  boundary conditions (BCs)
\begin{equation}\label{bsm-bc}
\psi_{\ell} \Big|_{\stackrel{\mu>0}{x=0}} =  \psi_{\ell}^{in+} \, , \quad
\psi_{\ell}\Big|_{\stackrel{\mu<0}{x=X}} =  \psi_{\ell}^{in-} \, ,
\end{equation}
where
$x$ is the spatial position;
$\mu$ is   the directional cosine of particle motion;
$\psi_{\ell}(x,\mu)$  is the conditional ensemble average  (CEA) of the angular flux provided  that
 $x$ is  in  material $\ell$. \linebreak
$p_{\ell}$ is the probability that the $\ell$-th material is found at a  position in the spatial domain.
For  homogeneous mixing statistics, the probability is given by
\begin{equation}
p_{\ell} = \frac{\lambda_{\ell}}{\lambda_1 + \lambda_2} \, .
\end{equation}
The  total ensemble average of the angular flux is defined by
\begin{equation}
\big< \psi \big>   \stackrel{\Delta}{=}  p_1 \psi_1 + p_2 \psi_2 \, .
\end{equation}
The CEA of the scalar flux and current for  material $\ell$ are given by
\begin{equation}
 \phi_{\ell} \stackrel{\Delta}{=} \int_{-1}^{1} \psi_{\ell} d \mu   \, , \quad
 J_{\ell} \stackrel{\Delta}{=} \int_{-1}^{1} \mu \psi_{\ell} d \mu \, .
\end{equation}

The problems with randomly mixed materials arise in various application, for example,
 in modeling  of inertial confinement fusion targets,
 nuclear fuel and shielding materials,
 clouds  in atmospheric sciences,
 tissue in radiation treatment planning \cite{pomraning-book}.
The source iteration scheme for the  linear binary stochastic transport equation
converges very slow \cite{palmer-tans-2000,miller-ttsp-2002}.
Various  synthetic acceleration methods
for this class of particle transport problems  have been developed \cite{adams-1991,prinja-fichtl-nse-2007,palmer-ane-2008,palmer-nse-2008,fichtl-jcp-2009}.
In this paper, we apply  nonlinear projection approach to derive
a multilevel iteration method \cite{Goldin-sbornik-82,dya-vyag-ttsp-1993}. It is defined by a system of equations consisting of the high-order transport equation and low-order  equations for various moments of the transport solution. The iteration method can be interpreted as a nonlinear multigrid
method over elements of the phase  space.
The proposed multilevel method can be applied for  radiation transport problems in random media  coupled with multiphysics equations \cite{miller-jqsrt-2001,dya-jcp-2019}.

The reminder of the paper is organized as follows.
 In Sec. \ref{sec:eqs}  the system of   low-order  equations  for    transport problems in BSM  is formulated.
 The iteration algorithm is described in Sec. \ref{sec:iter}.
 Numerical results are presented in Sec. \ref{sec:num}.
 We conclude with a brief discussion in Sec. \ref{sec:end}.

\section{\label{sec:eqs} Hierarchy of Equations of Multilevel Method for BSM}

The multilevel method  for solving the BSM transport problem \eqref{bsm-te} and \eqref{bsm-bc} is formulated by a  high-order equation for $\psi_{\ell}$ and low-order equations for  moments the  CEA of the angular flux  with exact closures.
 To derive the  hierarchy of equations of the method
 we  define (i)  projection operators   at each level,
  (ii)  moments of the transport solution    low-order equations are formulated for,
 and  (iii) prolongation  operators that enable to couple equations at different levels.

A specific feature of  the transport equation \eqref{bsm-te} for each material is that its right-hand side (RHS) contains  terms with absolute value of $\mu$ while the left-hand side (LHS) explicitly depends of $\mu$.
A natural way of forming a set of moment equations
based on projection of the high-order equation with this structure of  terms is
to define the low-order  equations for half-range moments
\begin{equation}
\phi_{\ell}^{\pm} \stackrel{\Delta}{=}   \pm\int_{0}^{\pm 1 } \psi_{\ell} d \mu
\end{equation}
that are the CEA  of the partial scalar fluxes.
The high-order equation \eqref{bsm-te}  is integrated over \linebreak
$-1 \le \mu \le 1$ with weights 1 and $\mu$ to obtain the low-order equations for $\phi_{\ell}^{\pm}$ defined by
\begin{subequations}  \label{loym}
\begin{multline} \label{loym0}
\frac{d}{dx} \big(p_{\ell} ( C_{\ell}^+\phi_{\ell}^+ +   C_{\ell}^-\phi_{\ell}^- )\big)
+ \Big( \sigma_{a,\ell}  + \frac{C_{\ell}^+}{\lambda_{\ell}} \Big)  p_{\ell}  \phi_{\ell}^+
+ \Big( \sigma_{a,\ell}  - \frac{C_{\ell}^-}{\lambda_{\ell}} \Big)   p_{\ell} \phi_{\ell}^- = \\
 \frac{p_{\ell'}}{\lambda_{\ell'}} \big( C_{\ell'}^+ \phi_{\ell'}^+ - C_{\ell'}^- \phi_{\ell'}^- \big) +
p_{\ell} q_{\ell} \, ,
\end{multline}
\begin{multline}\label{loym1}
\frac{d}{dx} \big(p_{\ell} ( E_{\ell}^+\phi_{\ell}^+ +   E_{\ell}^-\phi_{\ell}^- ) \big)
+ \Big( \sigma_{t,\ell}C_{\ell}^+  + \frac{E_{\ell}^+}{\lambda_{\ell}} \Big)   p_{\ell}\phi_{\ell}^+
+ \Big( \sigma_{t,\ell}C_{\ell}^-  - \frac{E_{\ell}^-}{\lambda_{\ell}} \Big)  p_{\ell} \phi_{\ell}^- = \\
 \frac{p_{\ell'}}{\lambda_{\ell'}} \big( E_{\ell'}^+ \phi_{\ell'}^+ - E_{\ell'}^- \phi_{\ell'}^- \big) \, ,
\end{multline}
\end{subequations}
where the low-order equations are closed by means of the  linear-fractional factors defined by
\begin{equation} \label{loym-factots}
C_{\ell}^{\pm} \stackrel{\Delta}{=} \frac{\int_0^{\pm 1} \mu \psi_{\ell'} d \mu}{\int_0^{\pm 1} \psi_{\ell} d \mu} \, ,
\quad
E_{\ell}^{\pm} \stackrel{\Delta}{=} \frac{\int_0^{\pm 1} \mu^2 \psi_{\ell} d \mu}{\int_0^{\pm 1} \psi_{\ell} d \mu} \, .
\end{equation}
The CEA of the material partial currents are defined by
\begin{equation} \label{j}
J^{\pm}_{\ell} \stackrel{\Delta}{=}   \pm\int_{0}^{\pm 1 } \mu \psi_{\ell} d \mu \, .
\end{equation}
They are cast in terms of the partial scalar fluxes as follows:
\begin{equation} \label{j-pm}
J^{\pm}_{\ell}  =    C^{\pm}_{\ell} \phi^{\pm}_{\ell} \, .
\end{equation}
The BCs  for Eqs. \eqref{loym} have the following form:
\begin{equation} \label{loym-bc}
\phi_{\ell}^+ \Big|_{x=0} =\phi_{\ell}^{in,+} \, ,
\quad
 \phi_{\ell}^- \Big|_{x=X} =\phi_{\ell}^{in,-} \, ,
\end{equation}
where
\begin{equation}
\phi_{\ell}^{in,\pm} \! \stackrel{\Delta}{=} \! \pm \int_0^{\pm 1}  \! \psi_{\ell}^{in,+} d \mu\,  .
\end{equation}
The low-order equations  \eqref{loym} use the same projection operators and unknowns
as the DP$_1$ equations \cite{yvon-1957,mertens-1954,bell,sanchez-2005}.
 They can be viewed as nonlinear  DP$_1$ equations with exact closures. This approach
 has been applied to formulate  Yvon-Mertens  (YM)  method for the linear Boltzmann equation  \cite{nikolai-1966}. Hereafter  these equations   are referred to  as
  the low-order YM (LOYM) equations.

The BSM transport equation \eqref{bsm-te} is now   projected  to formulate low-order equations  for the  total ensemble averaged  scalar flux $\big<\phi \big>$  and    current $\big< J \big>$, where
$\big<f \big> \stackrel{\Delta}{=} p_1 f_1 + p_2 f_2$
is the notation for the ensemble average.
The  equation \eqref{bsm-te} is summed over materials and then integrated over $-1 \le \mu \le 1$ with weights 1 and $\mu$ to obtain the
low-order quasidiffusion (LOQD) equations for $\big<\phi \big>$  and     $\big< J \big>$   given by \cite{gol'din-1964}
\begin{subequations} \label{loqd}
\begin{equation}
\frac{d\big< J \big>}{dx} + \big< \! \! \big< \sigma_a \big> \! \! \big> \big<\phi \big>  = \big< q \big> \, ,
\end{equation}
\begin{equation}
\frac{d}{dx} \Big(\big< \!\! \big<  E  \big> \! \! \big> \big< \phi \big> \Big) +\big< \!\! \big<  \sigma_t \big> \! \! \big>  \big< J \big>  +  \big< \! \! \big<  \eta  \big> \! \! \big>  \big< \phi \big>  = 0\, ,
\end{equation}
\end{subequations}
\begin{subequations} \label{qd-bc}
\begin{equation}
\big< J \big>\Big|_{x=0} \!\!  = \!\! \big< \! \! \big< C^- \big> \! \! \big> \Big( \big<\phi \big> - \big<\phi^{in,+} \big> \Big)\Big|_{x=0} \!\! +  \big<J^{in,+} \big> \, ,
\end{equation}
\begin{equation}
\big< J \big>\Big|_{x=X} \!\! =  \!\! \big< \! \! \big< C^+ \big> \! \! \big> \Big( \big<\phi \big> - \big<\phi^{in,-} \big> \Big)\Big|_{x=X}  \!\!+  \big<J^{in,-} \big> \, ,
\end{equation}
\end{subequations}
where the coefficients  are defined by means of the LOYM solution as follows:
\begin{equation}
\big< \!\! \big< \sigma_a \big> \! \! \big>  \stackrel{\Delta}{=} \frac{\sum_{\ell=1}^2  \sigma_{a, \ell} p_{\ell}   ( \phi_{\ell}^+  +  \phi_{\ell}^- ) }
{\sum_{\ell=1}^2 p_{\ell} ( \phi_{\ell}^+  +  \phi_{\ell}^-) } \, ,
\quad
\big< \!\! \big< \sigma_t \big> \! \! \big>  \stackrel{\Delta}{=} \frac{\sum_{\ell=1}^2  \sigma_{t, \ell} p_{\ell} |  C_{\ell}^+\phi_{\ell}^+  +  C_{\ell}^- \phi_{\ell}^-|  }
{\sum_{\ell=1}^2 p_{\ell} |  C_{\ell}^+\phi_{\ell}^+  +  C_{\ell}^- \phi_{\ell}^- | } \, ,
\end{equation}
\begin{equation}
\big< \!\! \big< \eta \big> \! \! \big>  \!  \stackrel{\Delta}{=}  \! \frac{\sum_{\ell=1}^2 \!   \big( \sigma_{t, \ell} \!  -  \!  \big< \! \! \big< \sigma_t \big> \! \! \big>  \big) p_{\ell}  (  C_{\ell}^+\phi_{\ell}^+  \!   +    \!  C_{\ell}^- \phi_{\ell}^-  )  }{\sum_{\ell=1}^2 p_{\ell}  ( \phi_{\ell}^+  +  \phi_{\ell}^-  ) } \, ,
\quad
\big< \!\! \big< E \big> \! \! \big>  \!   \stackrel{\Delta}{=} \!   \frac{  \sum_{\ell=1}^2  \!  p_{\ell}  ( E_{\ell}^+\phi_{\ell}^+ + E_{\ell}^-\phi_{\ell}^-  ) }
{\sum_{\ell=1}^2 p_{\ell}  ( \phi_{\ell}^+  +  \phi_{\ell}^-  )  } \, ,
\end{equation}
\begin{equation}
\big< \!\! \big< C^{\pm} \big> \! \! \big>\stackrel{\Delta}{=}  \frac{\sum_{\ell=1}^2     C_{\ell}^{\pm} p_{\ell} \phi_{\ell}^{\pm} }
{\sum_{\ell=1}^2 p_{\ell}  \phi_{\ell}^{\pm}  } \, ,
\quad
J_{\ell}^{in \pm} \stackrel{\Delta}{=} \int_{0}^{\pm  1} \mu \psi_{\ell}^{in \pm} d \mu \, .
\end{equation}

The next step in formulation of the multilevel system of low-order equations
is to define prolongation operators
that enable to couple the LOQD  and LOYM equations.
The term  $\frac{p_{\ell}}{\lambda_{\ell'}} ( C_{\ell}^+ \phi_{\ell}^+ - C_{\ell}^- \phi_{\ell}^-)$ is added to the LHS and RHS of Eq. \eqref{loym0} and
$\frac{p_{\ell}}{\lambda_{\ell'}} ( E_{\ell}^+ \phi_{\ell}^+ - E_{\ell}^- \phi_{\ell}^-)$
is added to  the LHS and RHS of Eq. \eqref{loym1} to get
\begin{subequations}  \label{loym-a}
\begin{multline} \label{loym0-a}
\frac{d}{dx} \big(p_{\ell} ( C_{\ell}^+\phi_{\ell}^+ +   C_{\ell}^-\phi_{\ell}^- )\big)
+ \bigg( \sigma_{a,\ell}  + C_{\ell}^+ \Big(\frac{1}{\lambda_{\ell}} + \frac{1}{\lambda_{\ell'}} \Big)  \bigg)  p_{\ell}  \phi_{\ell}^+
+ \bigg( \sigma_{a,\ell}  - C_{\ell}^- \Big(\frac{1}{\lambda_{\ell}} + \frac{1}{\lambda_{\ell'}} \Big)  \bigg)   p_{\ell} \phi_{\ell}^- = \\
 \frac{1}{\lambda_{\ell'}} \sum_{\ell''=1}^{2} p_{\ell''}\Big(  C_{\ell''}^+ \phi_{\ell''}^+ - C_{\ell''}^- \phi_{\ell''}^- \Big) +
p_{\ell} q_{\ell} \, ,
\end{multline}
\begin{multline}\label{loym1-a}
\frac{d}{dx} \big(p_{\ell} ( E_{\ell}^+\phi_{\ell}^+ +   E_{\ell}^-\phi_{\ell}^- ) \big)
+ \bigg( \sigma_{t,\ell}C_{\ell}^+  + E_{\ell}^+\Big(\frac{1}{\lambda_{\ell}} + \frac{1}{\lambda_{\ell'}} \Big)  \bigg)   p_{\ell}\phi_{\ell}^+
+ \bigg( \sigma_{t,\ell}C_{\ell}^-  - E_{\ell}^-\Big(\frac{1}{\lambda_{\ell}} + \frac{1}{\lambda_{\ell'}} \Big)  \bigg)  p_{\ell} \phi_{\ell}^- = \\
 \frac{1}{\lambda_{\ell'}} \sum_{\ell''=1}^{2} p_{\ell'}\Big( E_{\ell''}^+ \phi_{\ell''}^+ - E_{\ell''}^- \phi_{\ell''}^- \Big) \, .
\end{multline}
\end{subequations}
Taking into account Eq. \eqref{j-pm}, we obtain
\begin{equation} \label{c-phi}
 \sum_{\ell''=1}^{2} p_{\ell''}\Big(  C_{\ell''}^+ \phi_{\ell''}^+ - C_{\ell''}^- \phi_{\ell''}^- \Big) =
 \big< J^+ \big> -  \big< J^- \big> \, .
\end{equation}
We  define
\begin{equation} \label{<<E-mp>>}
\big< \!\! \big< E^{\pm} \big> \! \! \big> \stackrel{\Delta}{=}  \frac{\sum_{\ell=1}^2     E_{\ell}^{\pm} p_{\ell} \phi_{\ell}^{\pm} }
{\sum_{\ell=1}^2 p_{\ell}  \phi_{\ell}^{\pm}  } \,
\end{equation}
to get
\begin{equation} \label{e-phi}
 \sum_{\ell''=1}^{2} p_{\ell''}\Big(  E_{\ell''}^+ \phi_{\ell''}^+ - E_{\ell''}^- \phi_{\ell''}^- \Big) =
\big<\!\!\big<  E^+ \big> \!\!\big> \big<\phi^+ \big> - \big<\!\!\big<E^-\big>\!\! \big> \big< \phi^- \big> \, .
\end{equation}
Using Eqs. \eqref{c-phi} and \eqref{e-phi} in Eq. \eqref{loym-a}, we obtain the modified LOYM equations of the following form:
\begin{subequations}  \label{loym-mod}
\begin{multline} \label{loym0-mod}
\frac{d}{dx} \big(p_{\ell}( C_{\ell}^+\phi_{\ell}^+ +   C_{\ell}^-\phi_{\ell}^-) \big)
+ \bigg( \sigma_{a,\ell}  + C_{\ell}^+ \Big(\frac{1}{\lambda_{\ell}} + \frac{1}{\lambda_{\ell'}} \Big)  \bigg)   p_{\ell} \phi_{\ell}^+
+ \bigg( \sigma_{a,\ell}  -  C_{\ell}^- \Big(\frac{1}{\lambda_{\ell}} + \frac{1}{\lambda_{\ell'}} \Big) \bigg)   p_{\ell}\phi_{\ell}^- = \\
 \frac{1}{\lambda_{\ell'}} \Big(  \big< J^+ \big> - \big< J^- \big> \Big) +
p_{\ell} q_{\ell} \, ,
\end{multline}
\begin{multline}\label{loym1-mod}
\frac{d}{dx} \big(p_{\ell}( E_{\ell}^+\phi_{\ell}^+ +   E_{\ell}^-\phi_{\ell}^-) \big)
+ \bigg( \sigma_{t,\ell}C_{\ell}^+  + E_{\ell}^+  \Big(\frac{1}{\lambda_{\ell}} + \frac{1}{\lambda_{\ell'}} \Big) \bigg)  p_{\ell}  \phi_{\ell}^+
+ \bigg( \sigma_{t,\ell}C_{\ell}^-  -  E_{\ell}^-  \Big(\frac{1}{\lambda_{\ell}} + \frac{1}{\lambda_{\ell'}} \Big) \bigg)   p_{\ell}\phi_{\ell}^- = \\
 \frac{1}{\lambda_{\ell'}} \Big( \big<\!\!\big<  E^+ \big> \!\!\big> \big<\phi^+ \big> - \big<\!\!\big<E^-\big>\!\! \big> \big< \phi^- \big> \Big) \, .
\end{multline}
\end{subequations}
The prolongation operators that enable one to map
 $\big< \phi \big>$ and  $\big< J \big>$  to  $\big< \phi^{\pm} \big>$ and $\big< J^{\pm} \big>$
 are  defined  as follows \cite{dya-jsw-nse-2018}:
\begin{equation} \label{<phi+->}
\big< \phi^{\pm} \big> = \big< \! \! \big< \beta^{\pm} \big> \! \! \big> \Big(\big< J \big> -    \big< \! \! \big< \tilde C^{\pm} \big> \! \! \big>  \big<\phi \big>\Big) \, ,
\quad
\big< J^{\pm} \big> = \big< \! \! \big< \gamma^{\pm} \big> \! \! \big> \Big(\big< J \big> -    \big< \! \! \big< \tilde C^{\pm} \big> \! \! \big>  \big<\phi \big>\Big) \, ,
\end{equation}
where
\begin{equation}
\big< \!\! \big< \beta^{\pm} \big> \! \! \big> \stackrel{\Delta}{=}
 \frac{\sum_{\ell=1}^2       p_{\ell} \phi_{\ell}^{\pm} }
{\sum_{\ell=1}^2 ( C_{\ell}^{\pm} -  C_{\ell}^{\mp} )p_{\ell}  \phi_{\ell}^{\pm}  } \, ,
 \quad
 \big< \! \! \big< \gamma^{\pm} \big> \! \! \big> \stackrel{\Delta}{=}
 \frac{\sum_{\ell=1}^2     C_{\ell}^{\pm} p_{\ell} \phi_{\ell}^{\pm} }
{\sum_{\ell=1}^2 ( C_{\ell}^{\pm} -  C_{\ell}^{\mp} )p_{\ell}  \phi_{\ell}^{\pm}  } \, ,
\end{equation}
and
\begin{equation}
\big< \!\! \big< \tilde C^{\pm} \big> \! \! \big> \stackrel{\Delta}{=}  \frac{\sum_{\ell=1}^2     C_{\ell}^{\pm} p_{\ell} (\phi_{\ell}^{+} + \phi_{\ell}^{-}) }
{\sum_{\ell=1}^2 p_{\ell}  (\phi_{\ell}^{+} + \phi_{\ell}^{-}) }   \, .
\end{equation}

To formulate the high-order transport equation  of the multilevel method,
 we cast the RHS of Eq. \eqref{bsm-te} in terms of  $\phi_{\ell}$
to  obtain the    equation for $\psi_{\ell}$   given by
\begin{equation} \label{bsm-ho}
\mathcal{L}_{\ell} p_{\ell} \psi_{\ell}  =
\frac{1}{2}\sigma_{s, \ell} p_{\ell}   \phi_{\ell}
+ \frac{|\mu|}{\lambda_{\ell'} }p_{\ell'}  \psi_{\ell'}
+ \frac{1}{2}p_{\ell}  q_{\ell} \, ,  \ \ell' \ne \ell \, ,
\end{equation}
\begin{equation} \label{L-oper}
 \mathcal{L}_{\ell}  \stackrel{\Delta}{=}
\mu \frac{\partial      }{\partial  x}
+ \Big( \sigma_{t, \ell}  + \frac{|\mu|}{\lambda_{\ell}} \Big)   \, .
\end{equation}

\section{\label{sec:iter} Multilevel Iteration Scheme}

\begin{figure}[h!]
   \centering
    \includegraphics[width=0.25\textwidth]{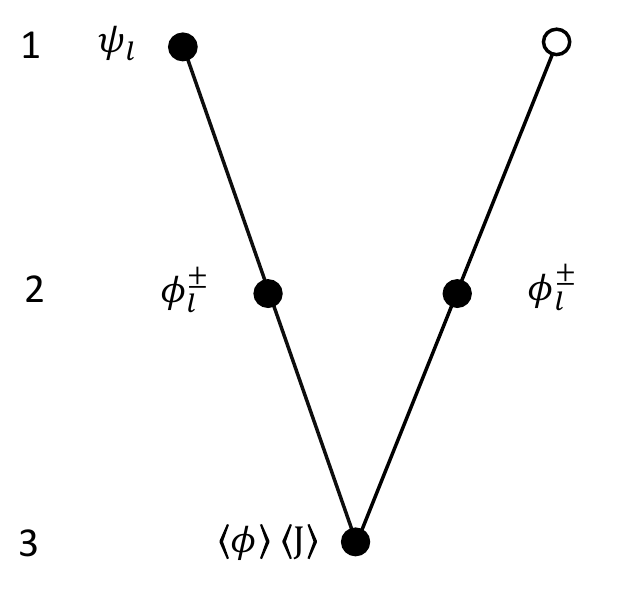}
  \caption{ \label{v-cycle} $V$-cycle}
 \end{figure}

The multilevel iteration method for  transport problems in BSM is defined
 by  several stages which involve solving the following  equations at different  levels:
\begin{itemize}
  \item the high-order transport equation \eqref{bsm-ho} for $\psi_{\ell}$,
      \item the LOYM equations \eqref{loym} and \eqref{loym-bc} for $\phi^{\pm}_{\ell}$,
  \item the LOQD equations \eqref{loqd} and \eqref{qd-bc} for $\big< \phi \big>$ and   $\big< J \big>$,
     \item the modified LOYM equations \eqref{loym-mod} and \eqref{loym-bc} for $\phi^{\pm}_{\ell}$.
\end{itemize}
This method is a nonlinear multigrid scheme over  elements of the phase space.
On each  transport  iteration, the hierarchy of equations is solved according to the $V$-cycle illustrated in Figure~\ref{v-cycle}.
The detailed  iteration scheme is  presented in Algorithm \ref{alg}, where $s$ is the index of transport iterations.

The  high-order  transport equation \eqref{bsm-ho} for  material $\ell$ depends on the CEA of the angular flux of another material $\ell'$. The LOYM   equations \eqref{loym} for  material $\ell$ also contain terms
  with the partial scalar  fluxes and factors of  material $\ell'$.
 In this paper,  fully iterative approach is studied.
 The multilevel method  executes   relaxation cycles at the level of  the material high-order equations using  Gauss-Seidel (GS) iterations. The  scattering term is defined by the material scalar fluxes obtained  from the low-order problems.  $n$ is the index of inner GS iterations.
 The number  of  GS cycles equals $n_{max}$.
 The same approach is applied to solve coupled material LOYM equations \eqref{loym}.
 At this level, only one GS iteration  over the material LOYM equations is performed.
 The modified  LOYM equations \eqref{loym-mod} for materials are decoupled
  due to the prolongation operator that defines the RHS. There are no inner iterations at this stage.

To discretize the high-order transport equation \eqref{bsm-ho} the linear discontinuous (LD)  finite element method is used.
The LOYM and LOQD equations  are discretized consistently with the LD scheme for Eq. \eqref{bsm-ho} \cite{dya-jsw-nse-2018}.

\begin{algorithm}[t!]
{
\DontPrintSemicolon
  $s=0$ \;
\While{$|| \big<\phi\big> ^{(s)}- \big<\phi\big> ^{(s-1)} || _{\infty}>
\epsilon || \big<\phi\big> ^{(s)}|| _{\infty}$ }{
 $s=s+1$ , \  $n=0$ , \ $\psi_{2}^{(0,s+1/3)} = \psi_{2}^{(s)}$\;
{{\bf Level  1}:  high-order transport   problem for materials}\;
\While{$n \le n_{max}$ }{
$n=n+1$\;
 $\mathcal{L}_1  p_1 \psi_{1}^{(n,s+1/3)}  =
|\mu|  \lambda_{2}^{-1}  p_{2}  \psi_{2}^{(n-1,s+1/3)}
+ 0.5\sigma_{s, 1} p_{1} \phi_{1}^{(s)}
+ 0.5p_{1}  q_{1} $\;
$\mathcal{L}_2  p_2 \psi_{2}^{(n,s+1/3)}  =  |\mu|  \lambda_{1}^{-1}  p_{1}  \psi_{1}^{(n,s+1/3)}
 \quad +0.5\sigma_{s, 2} p_{2} \phi_{2}^{(s)}
+ 0.5p_{2}  q_{2}$\;
}
$\psi_{\ell}^{(s+1/3)} \leftarrow \psi_{\ell}^{(n_{max},s+1/3)}$\;
Compute $C_{\ell}^{\pm (s+1/3)}$ and $E_{\ell}^{\pm (s+1/3)}$  ($\ell=1,2$), $\phi_{2}^{\pm(s+1/3)}$\;
 {{\bf Level 2}:  low-order   problem for materials}\;
 Perform one GS material iteration to solve LOYM eqs. \eqref{loym} for $\phi_{\ell}^{\pm(s+2/3)}$ with  $\phi_{2}^{\pm(s+1/3)}$ as an initial guess\;
Compute $\! \big< \!\! \big< \beta^{\pm} \big> \! \! \big>^{\! (s+2/3)}$,
$\big< \!\! \big< \gamma^{\pm} \big> \! \! \big>^{\! (s+2/3)}$,
$\big< \!\! \big< \tilde C^{\pm} \big> \! \! \big>^{\! (s+2/3)}$,
$\big< \! \!\big<   E^{\pm} \big> \! \! \big>^{\! (s+2/3)}$,
 $\big< \!\! \big< \sigma_a \big> \! \! \big>^{(s+2/3)}$,
 $\big< \!\! \big< \sigma_t \big> \! \! \big>^{(s+2/3)}$,
$\big< \!\! \big< \eta \big> \! \! \big>^{(s+2/3)}$,
 $\big< \!\! \big< E\big> \!\! \big>^{(s+2/3)}$,
$\big< \!\! \big< C^{\pm} \big> \! \! \big>^{(s+2/3)}$\;
{{\bf Level 3}:  low-order  problem for total ensemble average moments}\;
Solve LOQD eqs. \eqref{loqd} for $\big< \phi \big>^{(s+1)}$ and $\big< J \big>^{(s+1)}$ \;
Compute $\big< \phi^{\pm} \big>^{(s+2/3)}$ and $\big< J^{\pm}\big>^{(s+2/3)}$\;
{{\bf Level 2}:  low-order problem for materials}\;
Solve LOYM eqs. \eqref{loym-mod} for $\phi_{\ell}^{\pm(s+1)}$\;
Compute $\! \big< \!\! \big< \beta^{\pm} \big> \! \! \big>^{\! (s+1)}$,
$\big< \!\! \big< \gamma^{\pm} \big> \! \! \big>^{\! (s+1)}$,
$\big< \!\! \big< \tilde C^{\pm} \big> \! \! \big>^{\! (s+1)}$,
$\big< \! \!\big<   E^{\pm} \big> \! \! \big>^{\! (s+1)}$\;
Compute $\big< \phi^{\pm} \big>^{(s+1)}$ and $\big< J^{\pm}\big>^{(s+1)}$\;
  %
  {\bf Prolongation}\;
   $\phi_{\ell}^{(s+1)}   \leftarrow (\phi_{\ell}^{-(s+1)}+\phi_{\ell}^{+(s+1)})
 \big< \phi \big>^{(s+1)} \big/ \sum_{\ell'=1}^2 p_{\ell'} \big( \phi_{\ell'}^{- (s+1)} + \phi_{\ell'}^{+ (s+1)} \big)  $\;
$\tilde \phi^{\pm(s+1)}_{\ell}   \leftarrow \phi_{\ell}^{\pm(s+1)}
 \big< \phi^{\pm} \big>^{(s+1)} \big/ \sum_{\ell'=1}^2 p_{\ell'} \phi_{\ell'}^{\pm(s+1)}$\;
$\psi_{\ell}^{(s+1)} \leftarrow \pm \psi_{\ell}^{(s+1/3)} \tilde \phi_{\ell}^{\pm (s+1)} \big/
 \int_{0}^{\pm1}\psi_{\ell}^{(s+1/3)}  d \mu'$, \ $\mu \gtrless 0$  \;
}
}
\caption{\label{alg}  The  multilevel  iteration method  for the BSM transport equation }
\end{algorithm}

\section{\label{sec:num} Numerical Results}

To analyze   the  multilevel  iteration method we use four sets of numerical test problems
  \cite{mla-jqsrt-1989}.
The data for the tests are presented in Table \ref{tests}.
The problems are defined for a slab $0 \le x \le 100$  with vacuum boundaries.
The spatial mesh is uniform with   100   cells.
The double $S_4$ Gauss-Legendre  quadrature set is used.
The sets of tests  differ  from each other by  values of $\lambda_{\ell}\sigma_{t,\ell}$.
The problems in each set have materials with different  scattering ratios as well as total cross sections.

\begin{table}[t!]
\begin{center}
\caption{\label{tests} Data for Test Problems}
\begin{tabular}{|c|c|c|c|c|c|c|c|c|c|c| }
 \hline
Test  &  $\sigma_{t,1}$     &   $\frac{\sigma_{s,1}}{\sigma_{t,1}}$   &  $\lambda_1$   &  $q_1$       &   $\sigma_{t,2}$      &  $\frac{\sigma_{s,2}}{\sigma_{t,2}}$   &   $\lambda_2$   &  $q_1$   &   $\lambda_1 \sigma_{t,1}$  & $\lambda_2 \sigma_{t,2}$     \\ \hline \hline
A$_1$      & $\frac{10}{99}$     &  0     & $\frac{99}{100}$    & 1&   $\frac{100}{11}$    & 1     & $\frac{11}{100}$   & 1  &  0.1  & 1  \\ \hline
A$_2$      & $\frac{10}{99}$    &  1     & $\frac{99}{100}$    & 1&   $\frac{100}{11}$    & 0     & $\frac{11}{100}$   & 1  &  0.1  & 1 \\ \hline
A$_3$      & $\frac{10}{99}$    &  0.9  & $\frac{99}{100}$    & 1 &   $\frac{100}{11}$   & 0.9  & $\frac{11}{100}$   & 1  &  0.1   & 1   \\
     \hline \hline
B$_1$       & $\frac{10}{99}$   &  0     & $\frac{99}{10}$      & 1 &   $\frac{100}{11}$   & 1     & $\frac{11}{10}$    & 1  & 1  & 10   \\ \hline
B$_2$       & $\frac{10}{99}$   &  1     & $\frac{99}{10}$      & 1 &   $\frac{100}{11}$   & 0     & $\frac{11}{10}$    & 1  & 1  & 10    \\ \hline
B$_3$       & $\frac{10}{99}$   &  0.9   & $\frac{99}{10}$      & 1 &   $\frac{100}{11}$  & 0.9  & $\frac{11}{10}$     & 1  & 1  & 10   \\
  \hline \hline
C$_1$      & $\frac{2}{101}$   &  0     & $\frac{101}{20}$      & 1 &   $\frac{200}{101}$   & 1     & $\frac{101}{20}$    & 1  &  0.1  & 10  \\ \hline
C$_2$       & $\frac{2}{101}$   &  1     & $\frac{101}{20}$      & 1 &   $\frac{200}{101}$   & 0     & $\frac{101}{20}$    & 1  &  0.1  & 10   \\ \hline
C$_3$        & $\frac{2}{101}$   &  0.9  & $\frac{101}{20}$      & 1 &   $\frac{200}{101}$   & 0.9   & $\frac{101}{20}$    & 1  &  0.1  & 10   \\
  \hline    \hline
D$_1$       & $\frac{1}{99}$     &  0     & $\frac{99}{100}$      & 1 &   $\frac{10}{11}$      & 1     & $\frac{11}{1  0}$     & 1  & 0.01  & 1   \\ \hline
D$_2$      & $\frac{1}{99}$     &  1     & $\frac{99}{100}$      & 1 &   $\frac{10}{11}$      & 0     & $\frac{11}{1  0}$     & 1  & 0.01  & 1   \\ \hline
D$_3$      & $\frac{1}{99}$     &  0.9     & $\frac{99}{100}$      & 1 &   $\frac{10}{11}$      & 0.9     & $\frac{11}{1  0}$     & 1  & 0.01  & 1  \\ \hline
   \end{tabular}
\end{center}
\end{table}

\begin{table}[t!]
\begin{center}
\caption{\label{rho} Numerically Estimated Spectral Radii}
\begin{tabular}{|c|c|c|c|c|c|c|c|c|c|c|c|c| }
 \hline
$n_{max}$ & A$_1$  & A$_2$  &  A$_3$ & B$_1$ & B$_2$ &   B$_3$ & C$_1$ &  C$_2$ &  C$_3$ &  D$_1$ & D$_2$ & D$_3$  \\ \hline
1                & 0.4      &	0.41	   &	0.41   &	0.23	 &	0.12	 &	0.12      &	0.11	&	0.06    &	0.10     &  0.46     &	0.43    &	0.46 \\ \hline
2                & 0.16	 &	0.17	   &	0.16	 &	0.25	 &	0.12	 &	0.1        &	0.14	&	0.03    &	0.11     &  0.19     &	0.20    &	0.21 \\ \hline
   \end{tabular}
\end{center}
\end{table}

\begin{figure}[t]
	\centering \hspace*{-.5cm}
	\subfloat[Test A$_1$]{\includegraphics[width=.35\textwidth]{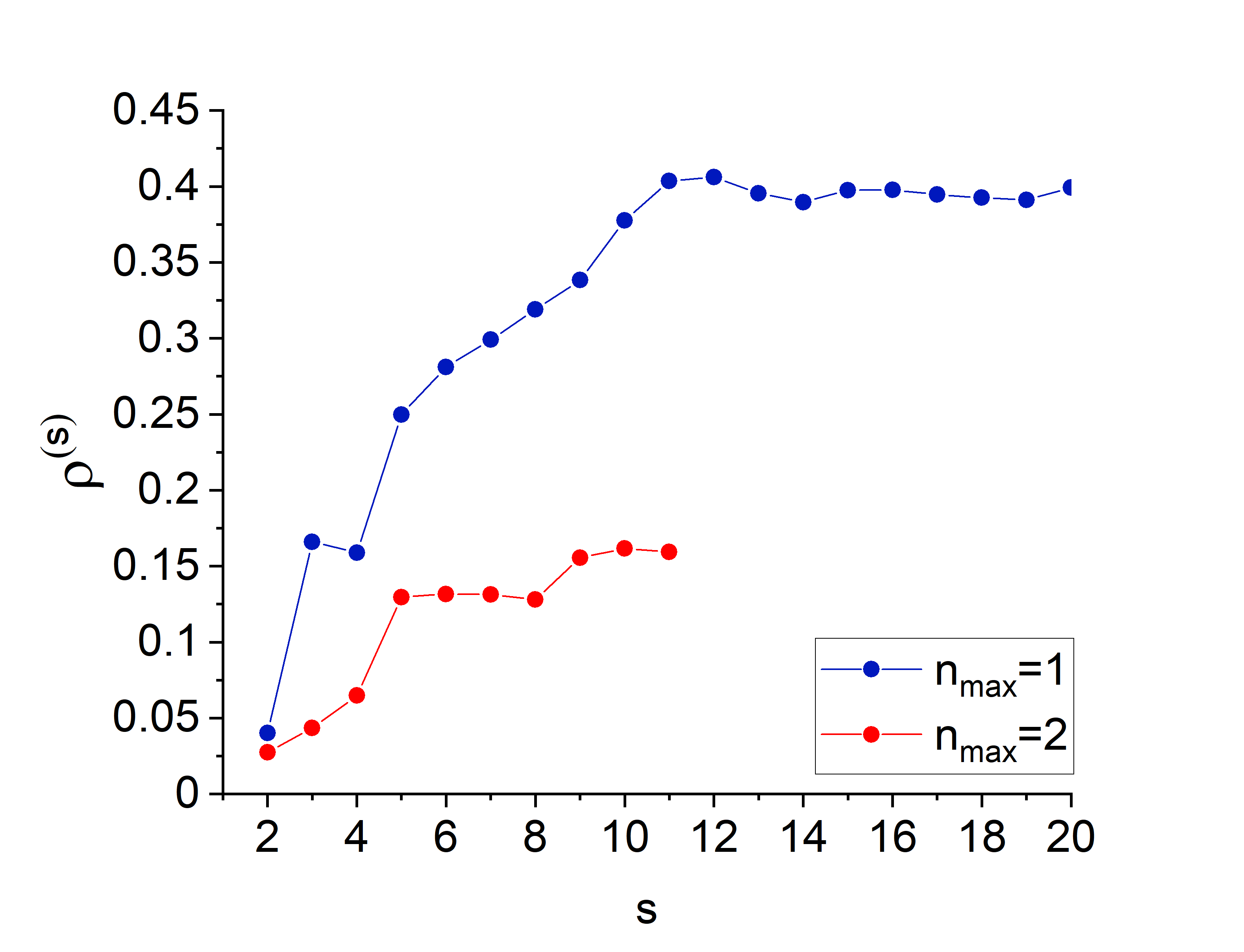}}
	\subfloat[Test A$_2$]{\includegraphics[width=.35\textwidth]{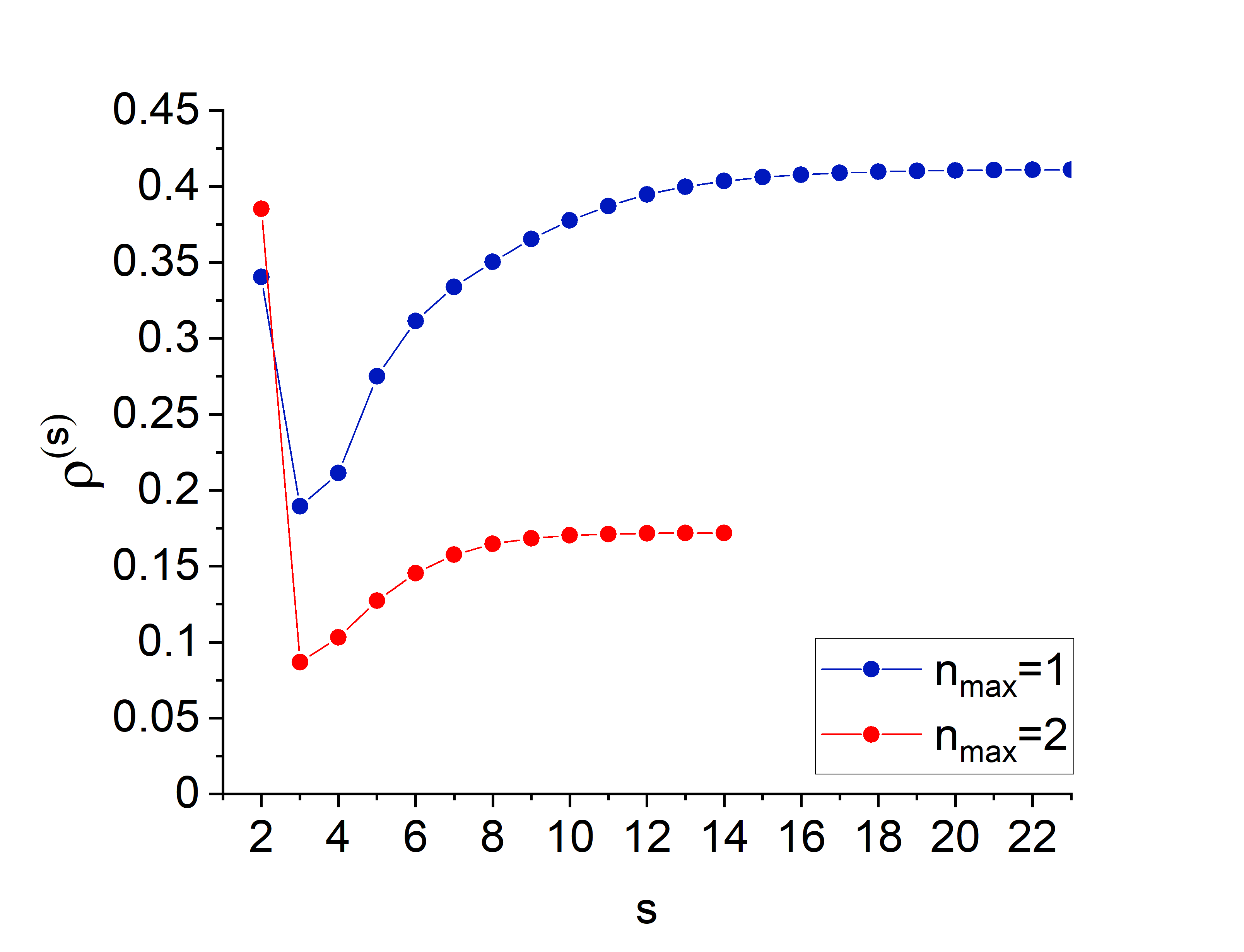}}
	\subfloat[Test A$_3$]{\includegraphics[width=.35\textwidth]{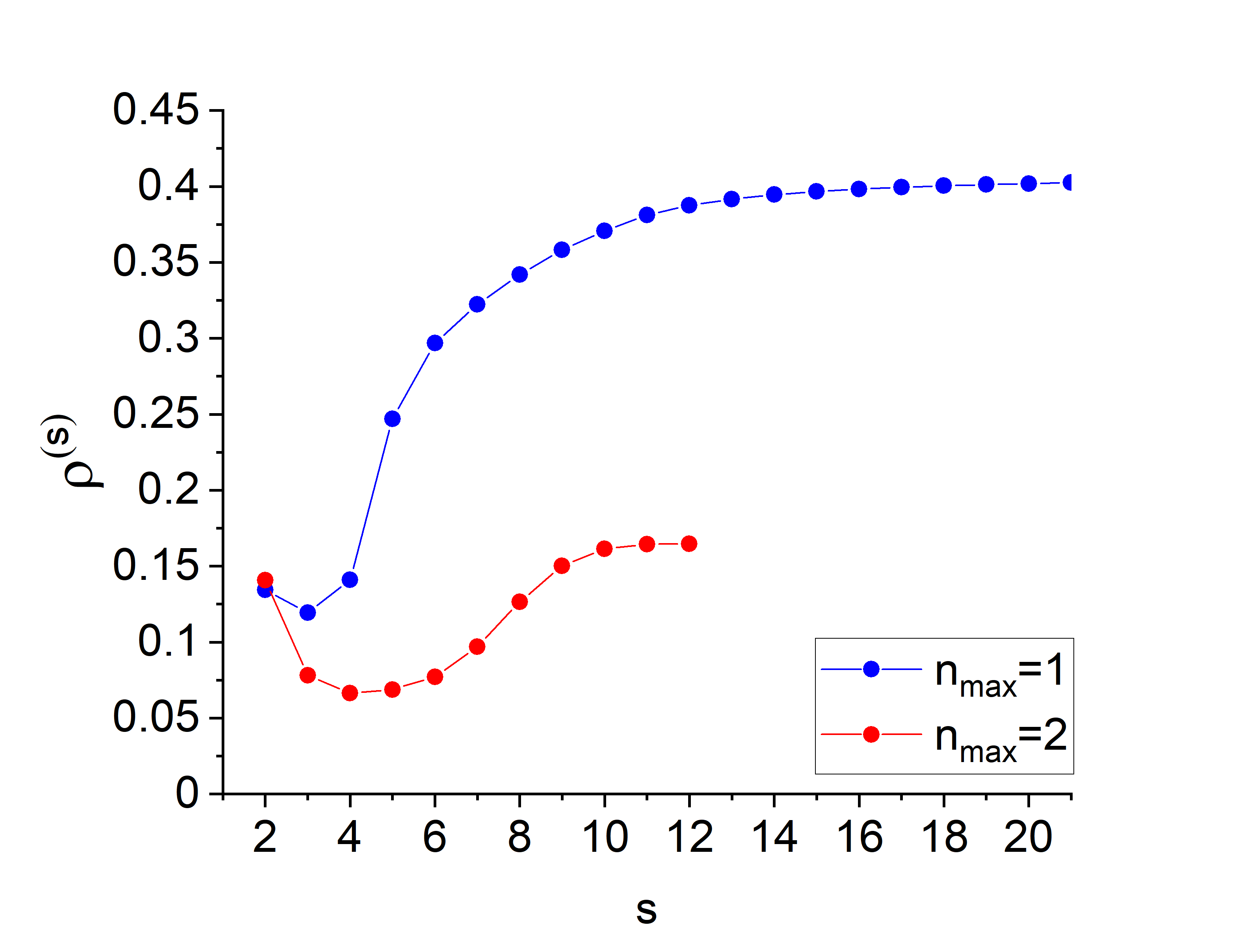}}
	\caption{\label{A-rho} Convergence rates $\rho^{(s)}$   in Tests A$_k$}

	\centering \hspace*{-.5cm}
	\subfloat[Test A$_1$]{\includegraphics[width=.35\textwidth]{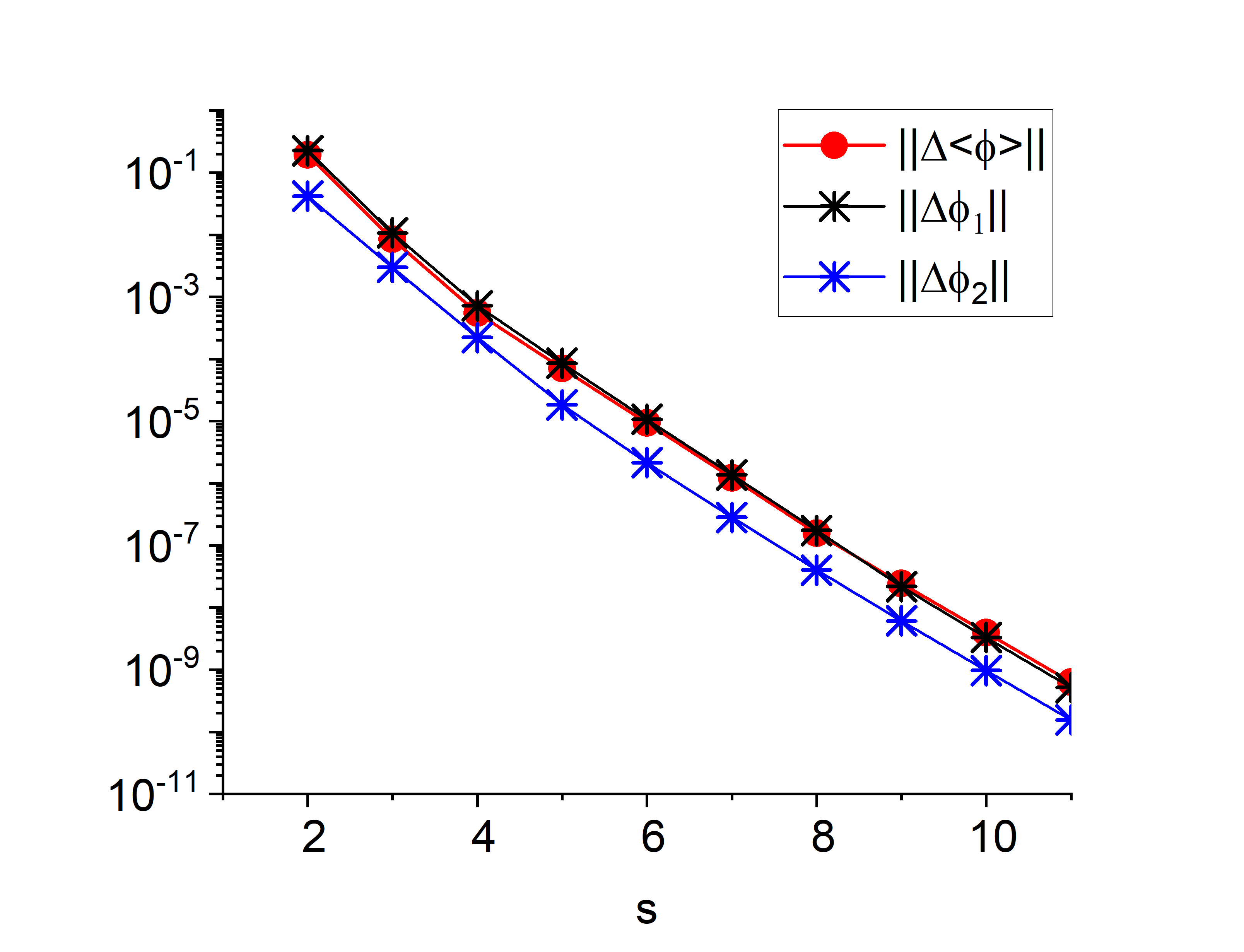}}
	\subfloat[Test A$_2$]{\includegraphics[width=.35\textwidth]{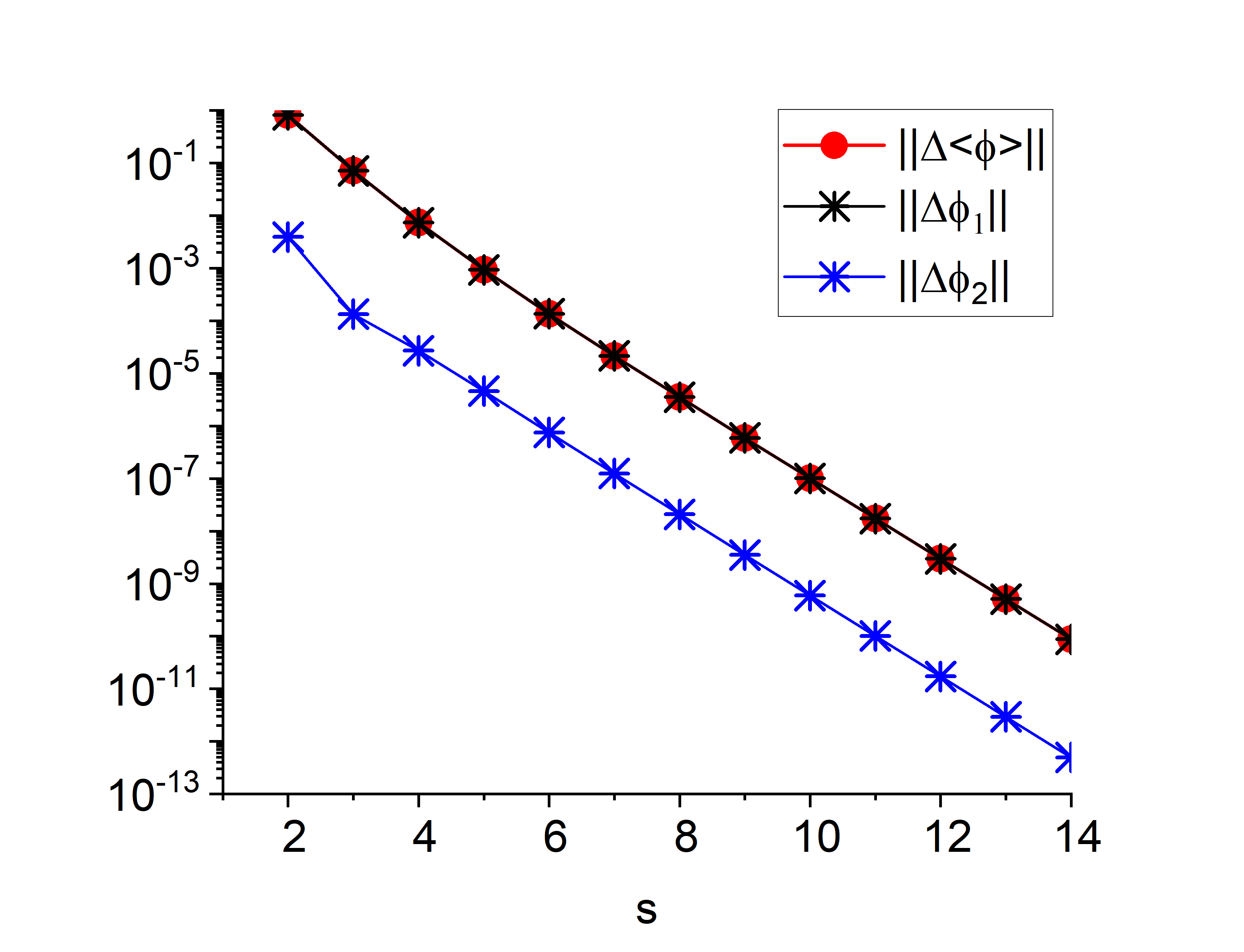}}
	\subfloat[Test A$_3$]{\includegraphics[width=.35\textwidth]{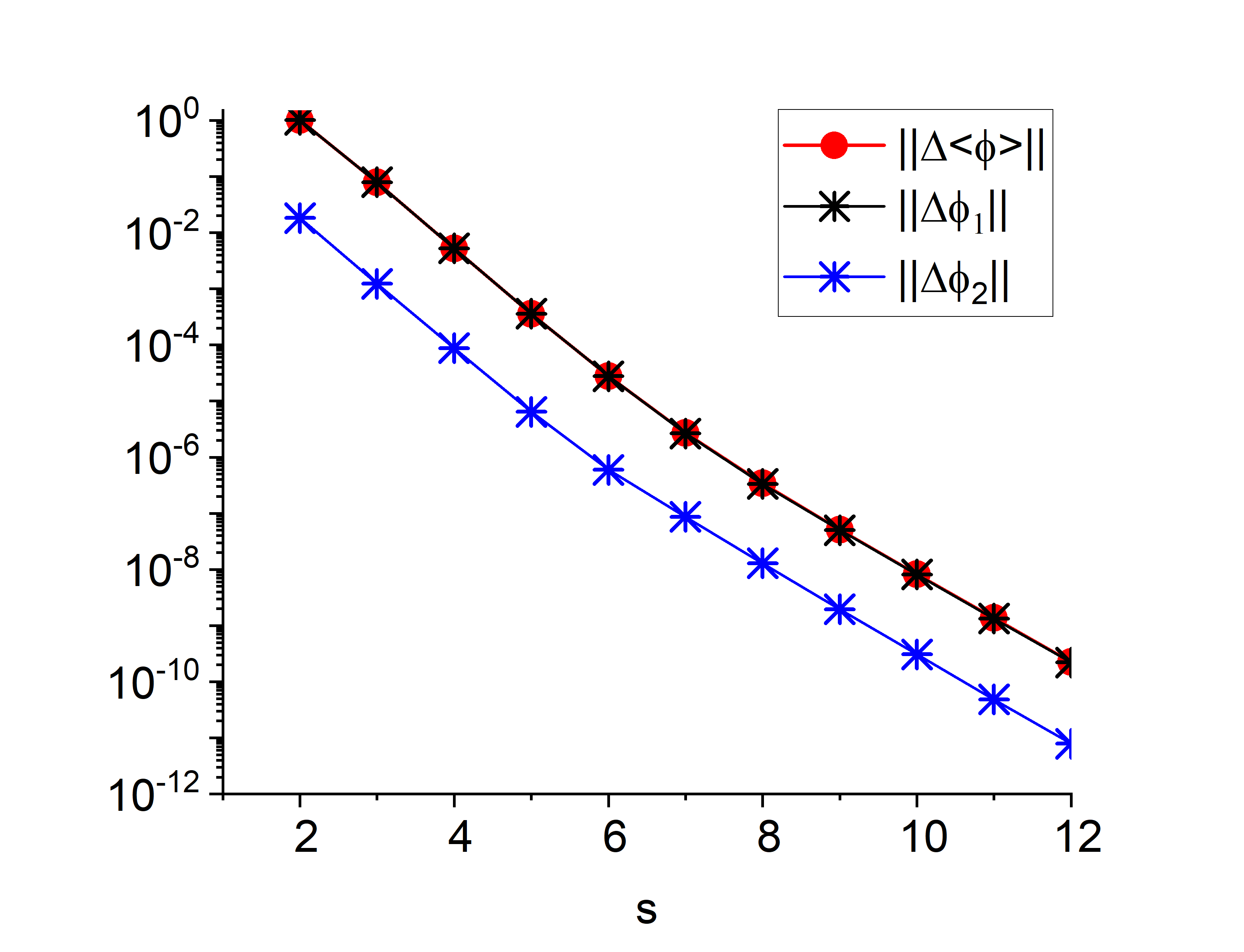}}
	\caption{\label{A-err} Convergence curves  of   $\big<\phi \big>^{(s)}$ and $\phi_{\ell}^{(s)}$ in Tests A$_k$, $n_{max}=2$}
 \end{figure}

\begin{figure}[t]

	\centering \hspace*{-.5cm}
	\subfloat[Test C$_1$]{\includegraphics[width=.35\textwidth]{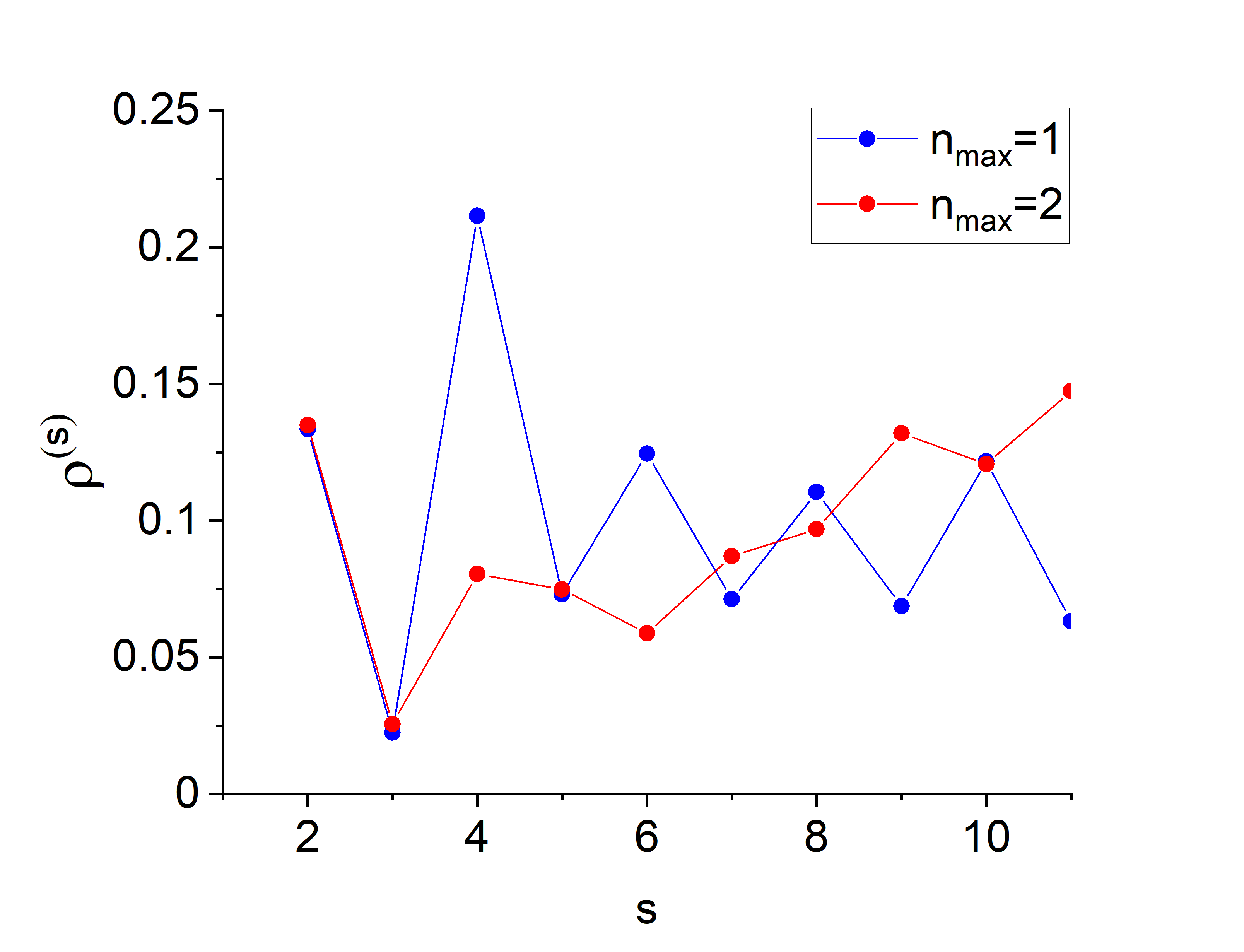}}
	\subfloat[Test C$_2$]{\includegraphics[width=.35\textwidth]{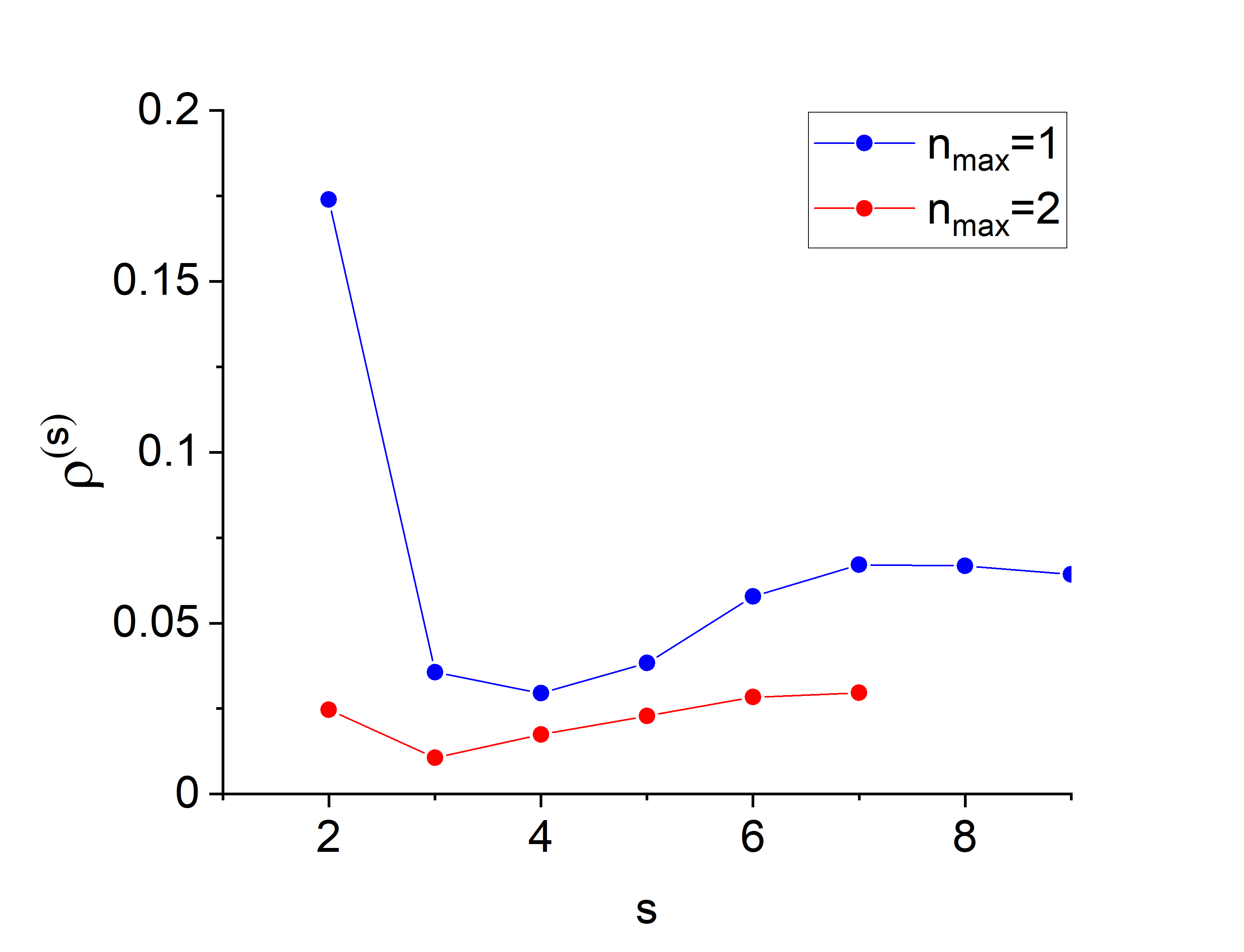}}
	\subfloat[Test C$_3$]{\includegraphics[width=.35\textwidth]{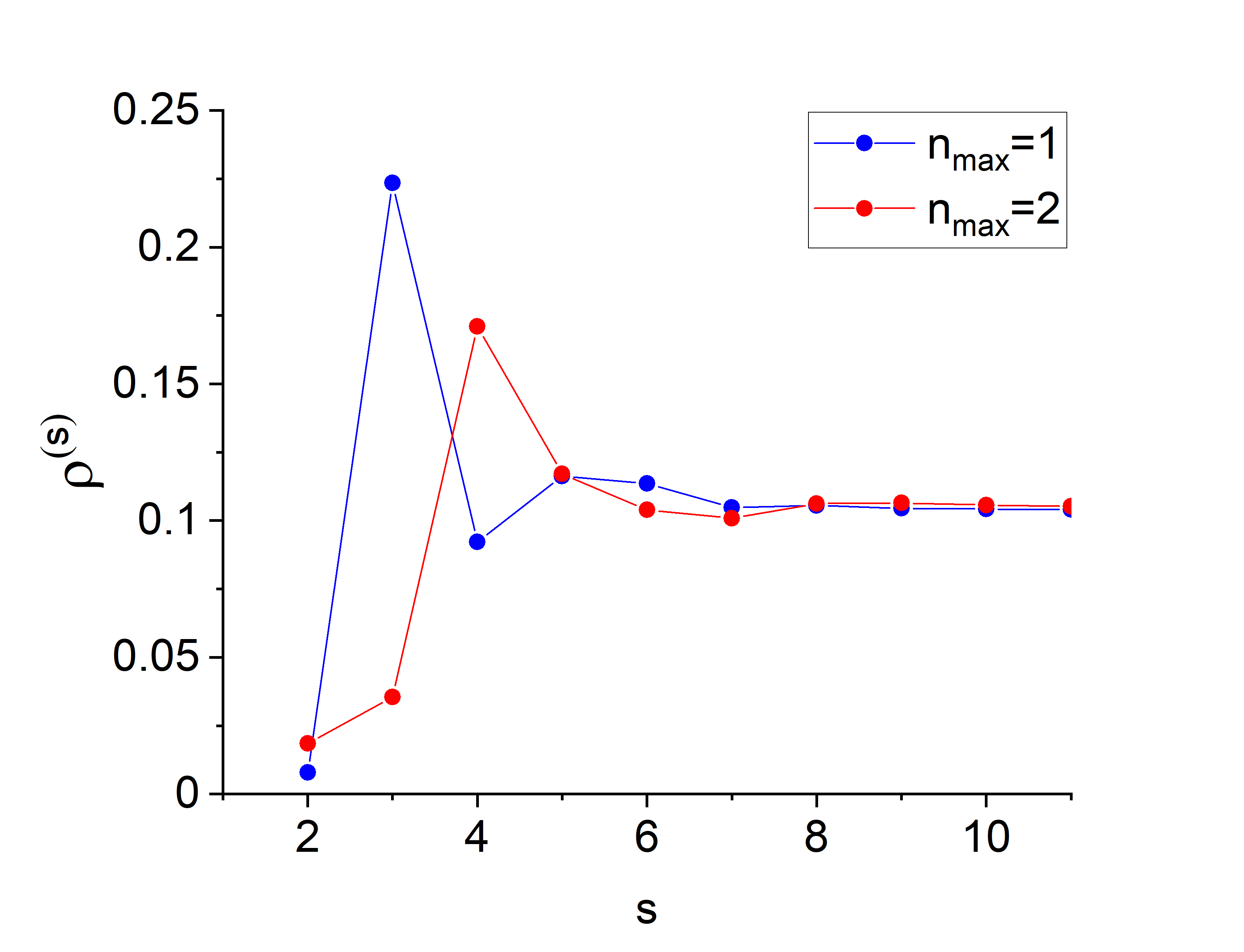}}
	\caption{\label{C-rho} Convergence rates $\rho^{(s)}$   in Tests C$_k$}

	\centering \hspace*{-.5cm}
	\subfloat[Test C$_1$]{\includegraphics[width=.35\textwidth]{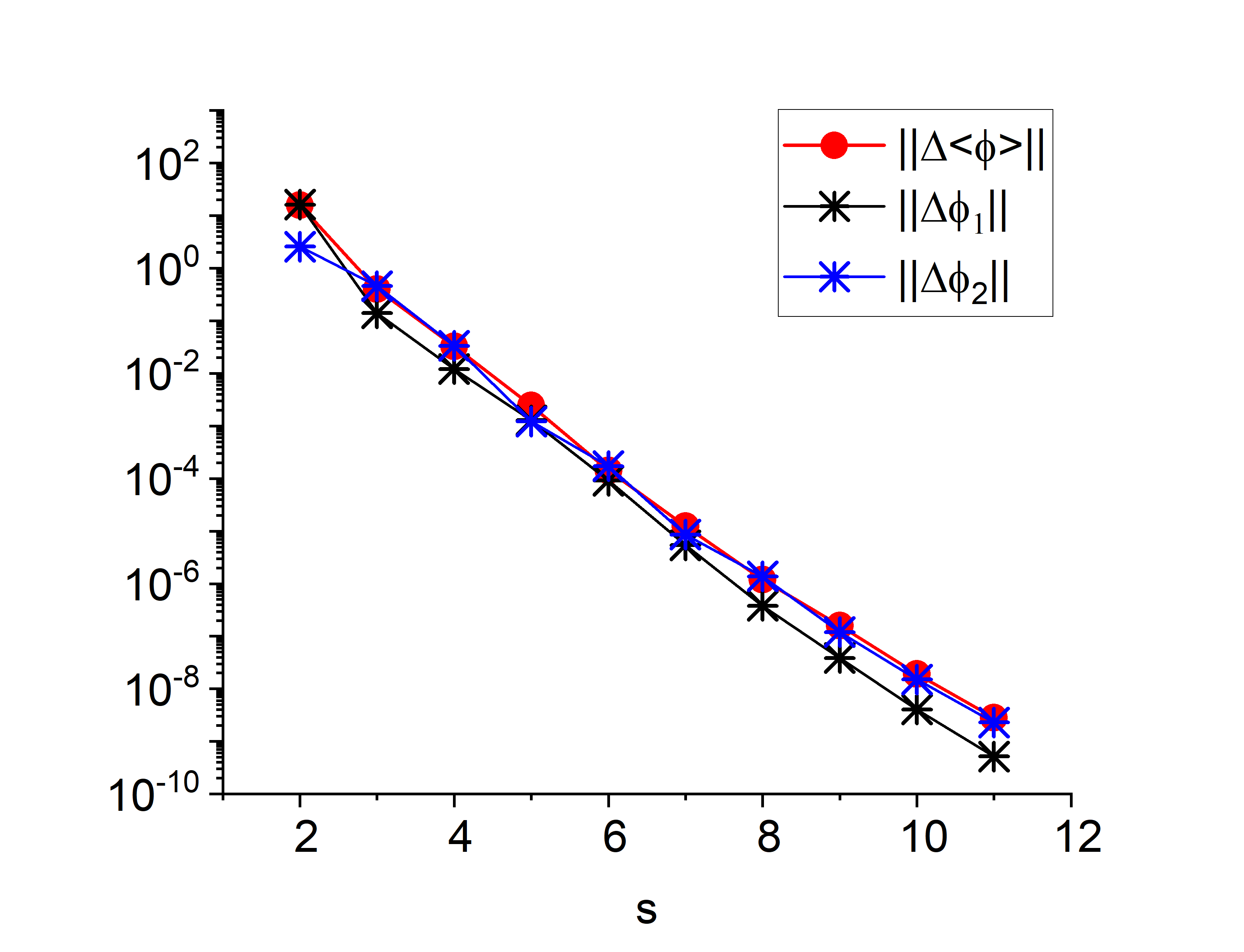}}
	\subfloat[Test C$_2$]{\includegraphics[width=.35\textwidth]{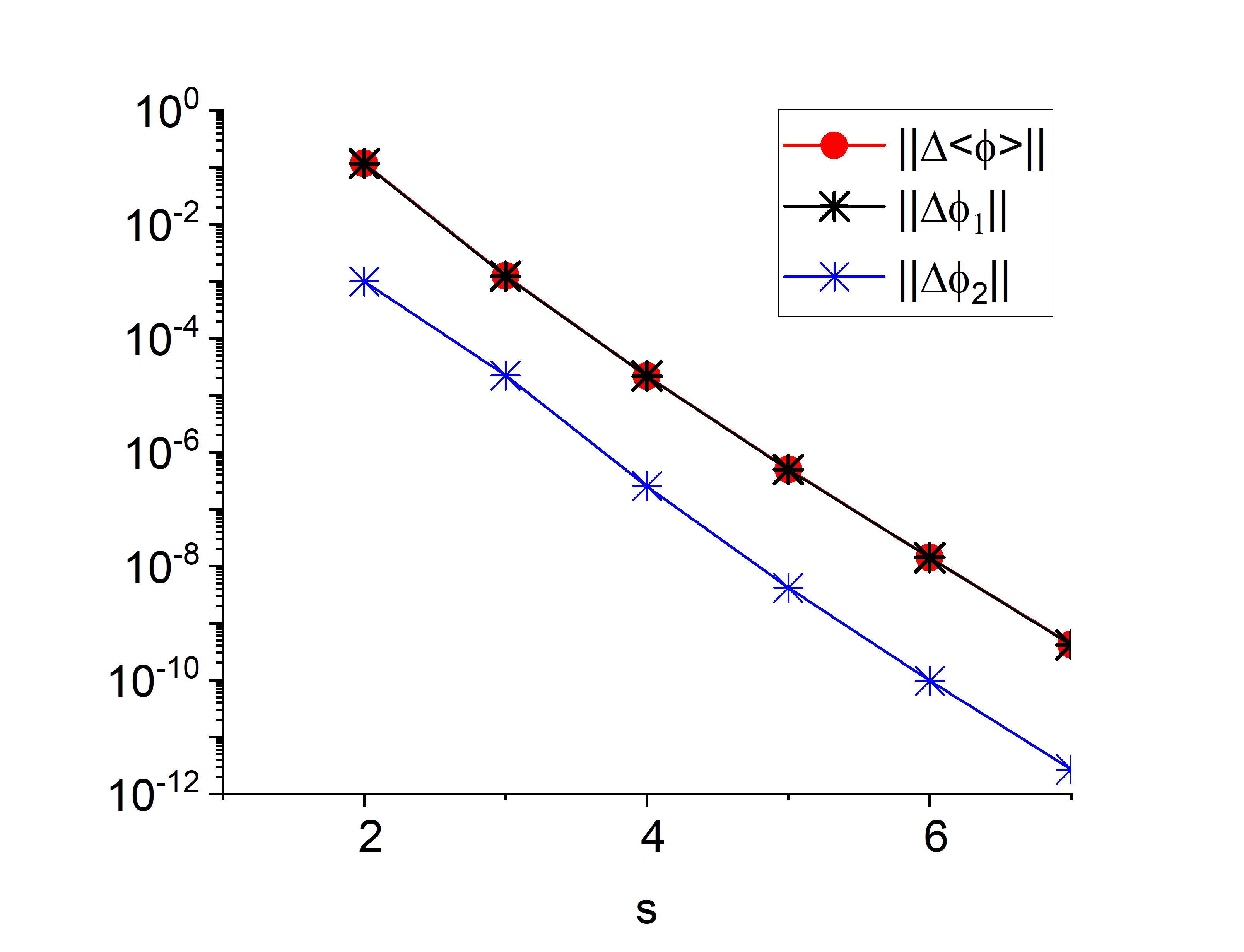}}
	\subfloat[Test C$_3$]{\includegraphics[width=.35\textwidth]{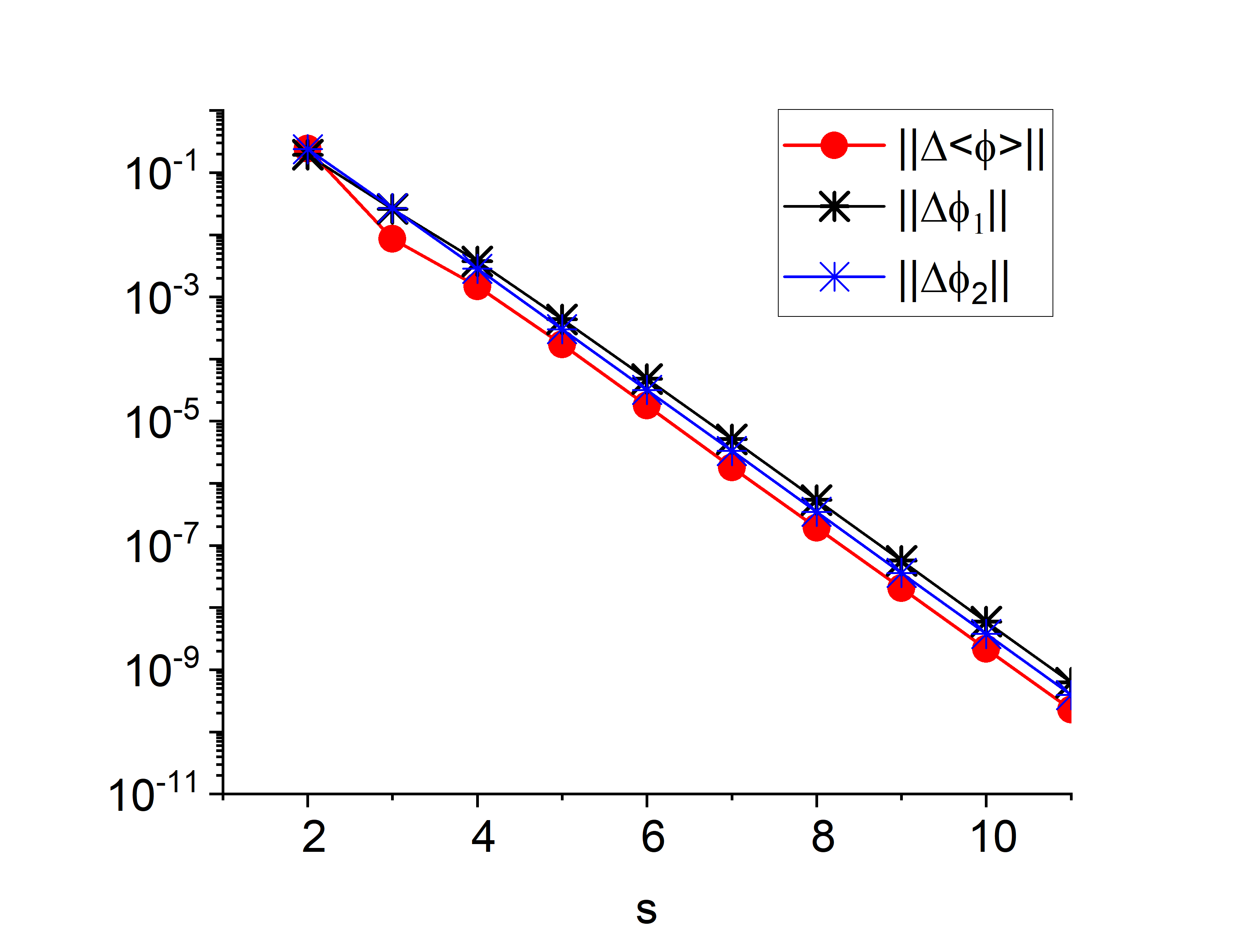}}
	\caption{\label{C-err} Convergence curves  of   $\big<\phi \big>^{(s)}$
and $\phi_{\ell}^{(s)}$ in Tests C$_k$, $n_{max}=2$}
\end{figure}

Table \ref{rho} presents numerically estimated spectral radii in every test.
The results for \linebreak
 algorithms with one and two  GS inner  iterations
 in the high-order transport problem \linebreak
 ($n_{max}=1,2$) are shown.
To demonstrate details of  iterative  behavior of the algorithms, Figures \ref{A-rho} and  \ref{C-rho}  show    convergence rates of the total ensemble average scalar flux
\begin{equation}
\rho^{(s)} = \frac{\big|\big| \Delta \big< \phi \big>^{(s)}\big|\big| _{\infty}}{\big|\big| \Delta \big< \phi \big>^{(s-1)}\big|\big| _{\infty}} \, ,
\quad
\Delta \big< \phi \big>^{(s)} = \big< \phi \big>^{(s)} -   \big< \phi \big>^{(s-1)} \,
\end{equation}
in Tests $A_k$ and $C_k$, respectively.
The algorithm with two relaxation cycles in the high-order problem ($n_{max}=2$)
converges fast in all tests.  In Tests $A_k$ and $D_k$, the spectral radii of this algorithm are smaller
than those of the algorithm with   one relaxation cycle ($n_{max}=1$).
In Tests $B_k$ and $C_k$, the algorithm
with one relaxation cycle  is rapidly converging.
Performance of the algorithms with 1  and 2 relaxation cycles are close to each other in these two sets of tests.

Figures \ref{A-err} and  \ref{C-err} present  convergence curves of
$||\Delta \big< \phi \big>^{(s)}||_{\infty}$ and
$|| \Delta  \phi_{\ell}^{(s)} || _{\infty}$
in Tests $A_k$ and $C_k$,  where
$\Delta \phi_{\ell}^{(s)} = \phi_{\ell}^{(s)} -   \phi_{\ell}^{(s-1)}$.
The CEA of  the material scalar fluxes
and   total ensemble average of the scalar flux converge with close rates.
The convergence behavior in other tests is similar.


\section{\label{sec:end} Conclusions}

A new  iteration method for  the transport equation in binary stochastic media has been presented.
It  is formulated by means of a hierarchy  of  equations consisting
of  the high-order transport equation for the CEA  of the angular flux and
the low-order equations for (i) the CEA of  the material partial scalar fluxes     and
 (ii) the  total ensemble average scalar flux and current.
The obtained results show that the multilevel method accelerates transport iterations
in the considered class of test problems.
Further analysis is needed to study efficiency of the iteration method.
One of items is  convergence in the atomic mix limit as
$\lambda_{\ell}\sigma_{t,\ell} \to 0$.

 There are various possible modifications that can be applied
to the multilevel method, for example, Anderson acceleration, nonlinear Krylov acceleration and advanced prolongation operators for the CEA of the angular flux.
The low-order equations of the multilevel method  correspond to different scales of the physical problem.
This method has a potential for  nonlinear radiation transport problems in random media  in which
the transport equation is coupled with multiphysics equations.

\bibliographystyle{elsarticle-num}
\bibliography{bsmedia-mc2021}

\end{document}